\documentclass[10pt]{article}
\usepackage[utf8]{inputenc}
\usepackage{amssymb}
\usepackage{amsmath}
\usepackage{amsthm}
\usepackage{tikz}
\usepackage{hyperref}
\usetikzlibrary{arrows,calc,decorations.pathreplacing,fadings,intersections}
\usetikzlibrary{external}
\newtheorem{theorem}{Theorem}
\newtheorem{lemma}[theorem]{Lemma}

\newtheorem{proposition}[theorem]{Proposition}
\newtheorem{definition}{Definition}

\newcommand{\vol}{\mathrm{vol}}

\title{%Relative Affine Surface Area \\
Floating bodies for ball-convex bodies}
\author{
Carsten Sch\"utt\footnote{Department of Mathematics, University of Kiel, Heinrich-Hecht-Platz 6, 24118 Kiel, Germany, 
{\tt schuett@math.uni-kiel.de }} \hskip 7mm
Elisabeth M. Werner \thanks{Supported by NSF grant DMS-2103482} \footnote{Department of Mathematics, Case Western Reserve University, 2145 Adalbert Road, Cleveland, OH 44106, USA,
{\tt elisabeth.werner@case.edu}} \hskip 7mm Diliya Yalikun \thanks{Supported by NSF grant DMS-2103482} \footnote{Department of Mathematics, Case Western Reserve University, 2145 Adalbert Road, Cleveland, OH 44106, USA, {\tt dxy259@case.edu@case.edu}} 
}
\date{}

\begin{document}

\maketitle

\begin{abstract}
We define floating bodies in the class of $n$-dimensional  ball-convex bodies. 
A right derivative of volume of these floating bodies leads to 
a  surface area measure for ball-convex bodies which
we call relative affine  surface  area.
We show that this  quantity is a rigid motion invariant, upper semi continuous valuation.
\end{abstract}

\section{Introduction}

A convex body  in $\mathbb R^{n}$ is a compact, convex subset of $\mathbb R^{n}$ with nonempty
interior. For a given convex body $C$,  we consider  the class of $C$-ball convex bodies $K$,  i.e., bodies 
that are either equal to the  empty set or to $\mathbb{R}^n$ or of the form
\begin{equation}\label{Cconvex}
K= 
\bigcap_{K \subset C+x} C+x.
\end{equation}
 This class was introduced and investigated  in \cite{LangiNazodiTalata}.
In this paper we will only consider  $C$-ball convex bodies $K$ that are non-empty and proper, i.e., $K \neq\emptyset$ is a proper subset of $\mathbb{R}^n$. 
\par
In particular, when  $C= R\, B^n_2$ is the $n$-dimensional Euclidean unit ball centered at $0$ with radius $R$,  we call $K$ $R$-ball convex and 
denote the class of such bodies by $\mathcal{K}_R$.
 This  class has been intensively studied, e.g., \cite{ArtsteinFlorentin, BezdekConnellyCsikos, BezdekNazodi, BezdekLangiNazodi, KabluchkoMarynychMolchanov, Pach}.
There are connections of this class to optimal transport, see \cite{ArtsteinSadovskyWyczensany}, to the Kneser-Poulsen conjecture (see \cite{Bezdek}) and to  isoperimetric problems, e.g., \cite{DrachTatarko}.
A classification of the isometries on ball-convex bodies  was given in  \cite{ArtsteinChorFlorentin1, ArtsteinChorFlorentin2}. 
\vskip 2mm
Two important closely related notions in convex geometry are the floating
body and the affine surface area of a convex body.
The floating body of a convex body is obtained by cutting off caps 
of volume less or equal to a fixed positive constant $\delta$.
Taking the right-derivative of the volume of the floating body gives rise to the affine surface area. 
This was established for all convex bodies in all dimensions 
in \cite{SW:1990}.
\par
The affine surface area $as(K)$ of a smooth convex body $K$ was introduced by Blaschke  \cite{Blaschke:1923}
as an integral over the boundary $\partial K$ of  $K$ of a power of the Gauss curvature $\kappa$, namely
\begin{equation}\label{asa}
as(K) = \int_{\partial K} \kappa(K, x)^\frac{1}{n+1} d\mu_K(x),
\end{equation}
where $\mu_K$ is the usual surface area measure on $\partial K$.
Due to its important properties, which make it an effective and powerful tool, it
is omnipresent in geometry, e.g., 
\cite{Lutwak:1996, LYZ:2000,  GHSW:2020, Hoehner:2022, LR:2010}).
There are numerous  applications for affine surface area, 
such as, 
$L_p$-extensions \cite{Lutwak:1996, SW:2004, SchuettWerner2023, TW:2019, TW:2023, WY:2008, Ye:2015, Ye:2016}, the approximation theory of convex bodies by polytopes 
\cite{BLW:2018, BGT, Boe1, 
%Gruber:1988,Gruber:1993, Ludwig:1999a, Boeroeczky:2000, Boeroeczky:2000a,
Reitz1, Reitzner, SW4}, affine curvature flows, information theory, and differential equations \cite{Andrews:1999, BesauWerner, TW1, TW2, TW4}.
\par
Extensions of both, floating body and affine surface area, to spherical and hyperbolic space were developed in \cite{BW:2016, BW:2018}.
\vskip 3mm
In this paper we consider the notion of $C$-ball floating body  for a $C$-ball-convex body.
In particular, taking the right-derivative of the volume when $C= R\, B^n_2$,  gives rise to  a  surface area measure for $R$-ball-convex bodies. 
As it  coincides with the  classical affine surface area when $R \to \infty$ and as  it  has similar properties 
we call it {\em relative affine surface  area}. 
The relative affine surface area is invariant under rigid motions, it  is a valuation and is upper semi-continuous.
In dimension $2$, this surface surface area measure has already appeared in approximation questions \cite{FodorKeveiVigh, FodorGrunfelder}.
\par
Further properties, as well as $L_p$-extensions, will be investigated in  forthcoming works. The authors believe that
both notions are of interest in its own right and will in particular be useful
for applications, such as approximation theory
of $C$-ball-convex bodies, extending the results that  have already been established in dimension $2$, e.g.,  \cite{FodorKeveiVigh, FodorGrunfelder}. 
\vskip 3mm
\noindent
{\bf Further notation.} 
The  closed Euclidean ball centered at $a$ with radius $r$ is $B^n_2(a,r)$. We write in short $B^n_2=B^n_2(0,1)$ and $S^{n-1}= \partial B^n_2$. The Euclidean norm on $\mathbb{R}^n$
is $\| \cdot\|$.
We denote by  $H\left(x, N\right)$ the hyperplane through $x$ orthogonal to the vector $N$. The line segment joining $x$ and $y$ is $[x,y]$.
Finally,  $c$, $d$ are absolute constants that may change from line to line.

\section{Results}

Let $K$ be a convex body in $\mathbb{R}^n$. For a point $x$ on the boundary $\partial K$ of $K$  we denote by $N_K(x)$ the outer unit normal vector of $K$ at $x$. 
By a theorem of Rademacher,  e.g., \cite{Rademacher}, $N_K(x)$ is unique a.e. on $\partial K$. The map $N: \partial K \to S^{n - 1}$ is called the spherical image map or Gauss map.  
Its  differential resp. generalized differential is called the Weingarten map. By a theorem of Alexandrov \cite{Alexandroff} and Busemann and Feller \cite{Buse-Feller}
these generalized derivatives exist a.e. on $\partial K$.
The eigenvalues of the Weingarten map are the principal curvatures $\kappa_i(K, x)$ of $K$ at $x$.
The (generalized) Gauss curvature is $\kappa(K, x) = \prod_{i=1}^{n-1} \kappa_i(K,x)$.
For more information and the details we refer to e.g., \cite{Gardner,  SchneiderBook}.
\vskip 4mm
We recall that for a convex body $K$ in $\mathbb{R}^n$, for $\delta \geq 0$, the (classical) convex floating body is defined as 
\begin{equation}\label{classicFB}
K_\delta= \bigcap _{\vol_n(H^- \cap K) \leq \delta} H^+.
\end{equation}
We now introduce  an analog to the classical floating body \cite{BaranyLarman1988, SW:1990}, the $C$-ball floating body.  
\vskip 2mm
\begin{definition}\label{def:R-float} 
Let $\delta \geq 0$. Let $K$ be a $C$-ball convex body.
 We define the $C$-ball floating body $K_\delta^C$ by 
\begin{equation}\label{L-float}
K_\delta^C = \bigcap_{\vol_n(K\setminus C+x) \leq \delta} C+x .
\end{equation}
In particular,  when $C=R\, B^n_2$, we  
 write in short $K_\delta^{R\, B^n_2}= K_\delta^R$ and call it the $R$-ball floating body of $K$.
\end{definition}
\par
Obviously, $K_\delta^C \subset K$ and $K_\delta^C$ is an $C$-ball convex body.
In this paper we will mostly treat the case $C=R\, B^n_2$.  Despite a slight abuse of notation, there should be 
no issue if we   
 write in short $K_\delta^{R\, B^n_2}= K_\delta^R$. 
\vskip 3mm
 In the definition of the $C$-convex floating body $K_\delta^C$ the convex body $C$ replaces the hyperplanes of the definition of the convex floating body (\ref{classicFB}). When $C=R\, B^n_2$ and $R \to \infty$, then $R \,B^n_2 \to H^+$ and we recover the usual convex floating body 
\cite{BaranyLarman1988, SW:1990}.
\vskip 3mm
Let  $g_{K}(x)=\vol_n(K\cap (x+K))$ be the  {\em covariogram operator}. 
We want to note that the $C$-ball floating body can be defined using an extension of the covariogram operator, the {\em cross variogram} operator
$$g_{K,C}(x)=\vol_n(K\cap (x+C)),$$
and then  
\begin{equation}\label{L-float2}
K_\delta^C = \bigcap \left\{ x +C: g_{K,C} (x) \geq \vol_n(K) -\delta\right\}.
\end{equation}
The covarigram is important in the recent development of {\em higher order convex bodies}, e.g.,  \cite{HaddadLangharstPuttermanRoysdonYe, LangharstXi}. It is also used to define the convolution bodies
\cite{Schm}, see also  \cite{W4},
\begin{equation*} \label{convolution}
    C(\delta,K)=\frac{1}{2} \left\{x\in\mathbb{R}^n:\vol_n(K\cap(x+K)) \geq \delta\right\} = \frac{1}{2} \left\{x\in\mathbb{R}^n: g_K(x) \geq \delta\right\}.
\end{equation*}
 I. Molchanov \cite{Molchanov} proposed an extension of the convolution bodies using the cross covariogram. 
It was this suggestion by I. Molchanov that inspired  definition (\ref{L-float2}). 
In dimension $2$, these are the $r$-spindle floating bodies \cite{FodorPapvari}.
\vskip 3mm
Note also that we  have for
all  affine maps $T$ on $\mathbb{R}^n$ with determinant $\mbox{det}(T)\neq 0$ 
that $K$ is $R$-ball convex if and only if $T(K)$ is $T(R\, B^n_2)$-ball convex
and that for all $\delta \geq 0$
\begin{equation}\label{Affine:map:floating:body}
(T(K))^{T(L)}_{\delta}=T\left(K^L_{\frac{\delta} {| det (T)|}} \right).
\end{equation}
\vskip 3mm
\noindent
The following theorem is the main result  of the  paper. We denote by 
$\mathcal{K}_R^+$ the set of all $R$-convex bodies in $\mathbb{R}^n$
such that for all $x \in \partial K$, $\kappa_i(K,x) >\frac{1}{R}$, for all $1 \leq i \leq n-1$. 
\vskip 2mm
\begin{theorem}\label{theorem:limit} 
Let $K \in \mathcal{K}_R^+$ and let $K_\delta^R$ be its $R$-ball floating body.  Then
$$
\lim_{\delta \to 0} \frac{\vol_n(K) - \vol_n(K_\delta^R)} {\delta^\frac{2}{n+1}} =  c_n \int_{\partial K} \prod _{i=1}^{n-1} \left( \kappa_i (K, x) -\frac{1}{R} \right)^\frac{1}{n+1} d \mu_K(x), 
$$
where $c_n=\frac{1}{2} \left( \frac{n+1}{\vol_{n-1}(B^{n-1}_{2})}\right)^\frac{2}{n+1}$.
\end{theorem}
\vskip 3mm
\noindent
Note that if $R \to \infty$, we recover the classical affine surface area (\ref{asa}).
This leads naturally to calling  the expression  of the above theorem {\em relative affine surface area}.
\begin{definition}\label{def:relative asa}
Let $K \in \mathcal{K}_R$.  Then the  relative affine surface area of $K$ is 
\begin{equation*}
as^R(K) = \int_{\partial K} \prod _{i=1}^{n-1} \left( \kappa_i (K, x) -\frac{1}{R} \right)^\frac{1}{n+1} d \mu_K(x)
\end{equation*}
\end{definition}
\vskip 3mm
\noindent
More generally,  when $K$ is an $L$-convex body, we define
\begin{equation*}
as^L(K) = \int_{\partial K} \prod _{i=1}^{n-1} \left( \kappa_i (K, x) - \kappa_i(L, N_L^{-1}(N_K(x)) \right)^\frac{1}{n+1} d \mu_K(x).
\end{equation*}
\vskip 2mm
\noindent
When $K$ is sufficiently smooth, $as^R$ can be written as
\begin{equation} \label{Sn}
as^R(K) = \int_{S^{n-1}} \prod _{i=1}^{n-1} \left( \frac{1}{r_i (K, u)} -\frac{1}{R} \right)^\frac{1}{n+1} f_K(u) d \sigma(u),
\end{equation}
and, more generally, 
\begin{equation*} \label{Sn}
as^L(K) = \int_{S^{n-1}} \prod _{i=1}^{n-1} \left( \frac{1}{r_i (K, u)} -\frac{1}{r_i(L,u)} \right)^\frac{1}{n+1} f_K(u) d \sigma(u),
\end{equation*}
where, for $u \in S^{n-1}$ such that $N_K(x)=u$,  $r_i (K, u)=(\kappa_i(K, N_K^{-1}(u))^{-1}$, and $f_K= \kappa^{-1}$ is the curvature function.
\vskip 3mm
\noindent
It cannot happen for an $R$-ball convex body that one of the $\kappa_i=0$. 
In fact,  for an $R$-convex body $\kappa_i \geq \frac{1}{R}$ for all $i$.
It can happen for an $R$-ball convex body that one or all of the $\kappa_i(x) =\frac{1}{R}$.
If one or all of the $\kappa_i(x)=\frac{1}{R}$ almost everywhere on $\partial K$, then $as^R(K)=0$.
For instance this is the case for {\em $R$-ball polyhedra},  i.e. the intersection of
finitely many $R$-balls. In particular $as^R(R B^n_2)=0$.
\vskip 3mm
\noindent
A crucial ingredient in the proof of Theorem \ref {theorem:limit} is the following proposition which  may be known but we could not find a proof.
It may also be useful in other contexts.
\vskip 3mm

\begin {proposition} \label{integral}
For $1 \leq i \leq n$, let $c_i >0$. Then
$$
\int_{S^{n-1}} \frac{d\sigma(\xi)}
{\left(\sum_{i=1}^{n} c_i \xi_i^2\right)^\frac{n}{2}} = \frac{2 \pi^\frac{n}{2}}{\Gamma(\frac{n}{2})} \left[\prod_{i=1}^{n} c_i ^\frac{1}{2}\right]^{-1} = \sigma (S^{n-1}) \left[\prod_{i=1}^{n} c_i ^\frac{1}{2}\right]^{-1}.
$$
\end{proposition}
\vskip3mm

We require in Proposition \ref{integral} that all $c_{i}>0$, $1\leq i\leq n$. If only one of the $c_{i}$
equals $0$ the integral equals $\infty$. Indeed, suppose that $c_{n}=0$, while the others are strictly greater than $0$. Then on a set of volume $d(n)\,\varepsilon^{n-1}$ the integrand $\left(\sum_{i=1}^{n} c_i \xi_i^2\right)^\frac{n}{2}$ is greater than $\frac{1}{\varepsilon^{n}}\left(\sum_{i=1}^{n} c_i\right)^\frac{n}{2}$. Thus the integral is unbounded.
\vskip 3mm
It is obvious that relative affine surface area is invariant under rigid motions. The next proposition states some more properties. Further properties will be in a forthcoming paper.
We will need the definition of the Hausdorff distance $d$ of convex sets $K$ and $L$,
$$
d(K, L) = \inf\{\lambda >0: K \subset L+\lambda B^n_2, L \subset K +\lambda B^n_2\}.
$$
\par
\begin{proposition}\label{rel-asa-properties} 
Let $K$ be an $R$-ball convex body. 
\vskip 1mm
\noindent
(i) Let $a \in \mathbb{R}$, $a >0$.  Then $as^R$ is homogeneous of degree $n \frac{n-1}{n+1}$, 
\begin{equation*}
as^{a R} (a K) = a^{n \frac{n-1}{n+1}} \, as^R(K).
\end{equation*}
\vskip 1mm
\noindent
(ii) The relative affine surface area $as^R$ is a valuation, i.e. for $R$-ball convex bodies $K$ and $L$ such that $K\cup L$ is again $R$-ball convex, 
$$
as^R(K\cup L) + as^R(K\cap L) = as^R(K) +as^R(L).
$$
\vskip 1mm
\noindent
(iii) We have
\begin{equation*}
as^R(K) \leq n \, \vol_n(B^n_2)^\frac{2}{n+1} \vol_n(K)^\frac{n-1}{n+1},
\end{equation*}
with equality if and only if $R=\infty$ and $K$ is an ellipsoid.
\vskip 1mm
\noindent
(iv) The relative affine surface area $as^R$ is upper semi continuous: If $K_l$, $l \in \mathbb{N}$,  is a sequence of $R$-ball convex 
bodies that converge in the Hausdorff distance to the $R$-ball convex body $K$, then
$$
\lim \sup_{l \in \mathbb{N} }\, as^R(K_l) \leq as^R(K).
$$
\end{proposition}
\vskip 2mm
Relative affine surface area is not continuous: Every $R$-ball convex body with strictly positive relative affine surface area 
can be approximated by a sequence of $R$-ball polyhedra.

\vskip 4mm
\section{Proofs}
\vskip 3mm
\subsection{Proof of Proposition \ref{integral}}
\vskip 2mm
We state and prove the proposition in the form we will use it, namely for  $1 \leq i \leq n-1$, let  $c_i >0$.
We want to compute 
\begin{equation*}
I=\int_{S^{n-2}}\frac{d\xi}{(\sum_{i=1}^{n-1}c_i\xi_i^2)^{(n-1)/2}}.
\end{equation*}
We pass to polar coordinates.
Let
$
\operatorname{Pol}_{n-1}:[0,\pi]\times\cdots\times[0,\pi]\times
[0,2\pi]
\rightarrow\mathbb R^{n-1}
$
be given by
$$
\operatorname{Pol}_{n-1}(\phi_{1},\dots,\phi_{n-2})
=(\xi_{1}^{}(\phi_{1},\dots,\phi_{n-2}),\xi_{2}^{}(\phi_{1},\dots,\phi_{n-2}),
\dots,\xi_{n-1}^{}(\phi_{1},\dots,\phi_{n-2}))
$$
and
\begin{eqnarray}\label{Polarkoordinaten-20}
\xi_1^{}(\phi_{1},\dots,\phi_{n-1})  & =&\cos \phi_{1}   \nonumber  \\
\xi_{k}^{}(\phi_{1},\dots,\phi_{n-1}) & =&\left(\prod_{i=1}^{k-1}\sin\phi_{i}\right)\cos\phi_{k}
\hskip 20mm
k=2,3,\dots,n-2     \\
\xi_{n-1}^{}(\phi_{1},\dots,\phi_{n-1}) & =&\prod_{i=1}^{n-2}\sin\phi_{i}
\nonumber
\end{eqnarray}
The factor due to the change of integration is
\begin{eqnarray}\label{Polarkoordinaten-5}
\prod_{k=1}^{n-2}\sin^{n-k-2}\phi_{k}.
\end{eqnarray}
We have
\begin{eqnarray}\label{Integral1}
I&=&\int_{S^{n-2}}\frac{d\xi}{(\sum_{i=1}^{n-1}c_i\xi_i^2)^{(n-1)/2}}
\\
&=&\int_0^{2\pi}\int_0^\pi\cdots \int_{0}^{\pi}
\frac{\prod_{k=1}^{n-2}\sin^{n-k-2}\phi_{k}\hskip3mm d\phi_{1}\cdots d\phi_{n-2}}
{\left(c_{1}\cos^{2} \phi_{1}+\sum_{k=2}^{n-2}c_{k}
\left(\prod_{i=1}^{k-1}\sin\phi_{i}\right)^{2}\cos^{2}\phi_{k}
+c_{n-1}(\prod_{i=1}^{n-2}\sin\phi_{i})^{2}\right)^{\frac{n-1}{2}}}
\nonumber
\end{eqnarray}
Now we successively integrate, starting with the integration in the variable $\phi_{1}$.
This leads to having  to consider successively integrals of the form
$$
J=\int_0^\pi\frac{\sin^{m}(x)}{(a\cos^2 x+b\sin^2 x)^{(m+2)/2}}dx, \hskip 4mm 0 \leq m \leq (n+2)/2.
$$
 Note that $\sin^m(\pi-x)=\sin^{m}x$ and $\cos^2(x)=\cos^2(\pi-x)$, so $J$ is symmetric about $\pi/2$.  Thus it is enough to consider
 $$
 J'=\frac{J}{2} =\int_0^{\pi/2}\frac{\sin^{m}(x)}{(a\cos^2x+b\sin^2x)^{(m+2)/2}}dx.
 $$ 
 We make a change of variables, $u=\frac{b}{a} \, \tan^2x$ and the integral becomes      
$$
  J'=\frac{1}{2 \,a^\frac{1}{2}\, b^{\frac{m+1}{2}}}\int_0^\infty u^{\frac{m-1}{2}}(1+u)^{-\frac{m+2}{2}}du 
 $$
 We use $u=\frac{1}{t}-1$, i.e. $t=\frac{1}{1+u}$ and $du=-\frac{1}{t^{2}}dt$
 \begin{eqnarray*}
 &&\frac{1}{2 \,a^\frac{1}{2}\, b^{\frac{m+1}{2}}}\int_0^{1} \left(\frac{1}{t}-1\right)^{\frac{m-1}{2}}t^{\frac{m+2}{2}}\frac{1}{t^{2}}dt
 =\frac{1}{2 \,a^\frac{1}{2}\, b^{\frac{m+1}{2}}} \int_0^{1}\left(1-t\right)^{\frac{m-1}{2}}t^{-\frac{1}{2}}dt
 \\
 &&=\frac{1}{2 \,a^\frac{1}{2}\, b^{\frac{m+1}{2}}}B\left(\frac{1}{2},\frac{m+1}{2}\right)
 =\frac{1}{2 \,a^\frac{1}{2}\, b^{\frac{m+1}{2}}}
 \frac{\Gamma(\frac{1}{2})\Gamma(\frac{m+1}{2})}{\Gamma(\frac{m+2}{2})}
  \end{eqnarray*}
and thus 
\begin{equation}\label{pre-integral}
J= \frac{\pi^\frac{1}{2}}{a^\frac{1}{2}\, b^{\frac{m+1}{2}}} 
  \frac{\Gamma(\frac{m+1}{2})}{\Gamma(\frac{m}{2}+1)}.
\end{equation}
Now we use (\ref{pre-integral}) with $m=n-3$ and integrate the integral \eqref{Integral1}
with respect to $\phi_{1}$
\begin{eqnarray*}
&&\int_0^{2\pi}\int_0^\pi\cdots \int_{0}^{\pi}
\frac{\prod_{k=1}^{n-2}\sin^{n-k-2}\phi_{k}\hskip3mm d\phi_{1}\cdots d\phi_{n-2}}
{\left(c_{1}\cos^{2} \phi_{1}+\sin^{2}\phi_{1}\left(\sum_{k=2}^{n-2}c_{k}
\left(\prod_{i=2}^{k-1}\sin\phi_{i}\right)^{2}\cos^{2}\phi_{k}
+c_{n-1}(\prod_{i=2}^{n-2}\sin\phi_{i})^{2}\right)\right)^{\frac{n-1}{2}}}
\\
&&=\frac{\sqrt{\pi}}{\sqrt{c_{1}}}\frac{\Gamma(\frac{n-2}{2})}{\Gamma(\frac{n-1}{2})}\int_0^{2\pi}\int_0^\pi\cdots \int_{0}^{\pi}
\frac{\prod_{k=2}^{n-2}\sin^{n-k-2}\phi_{k}\hskip3mm d\phi_{2}\cdots d\phi_{n-2}}
{\left(\sum_{k=2}^{n-2}c_{k}
\left(\prod_{i=2}^{k-1}\sin\phi_{i}\right)^{2}\cos^{2}\phi_{k}
+c_{n-1}\left(\prod_{i=2}^{n-2}\sin\phi_{i}\right)^{2}\right)^{\frac{n-2}{2}}}
\end{eqnarray*}
After having integrated with respect to $\phi_{1},\dots,\phi_{n-3}$, we arrive at
$$
I=\frac{\sqrt{\pi}}{\sqrt{c_1}}\frac{\Gamma(\frac{n-2}{2})}{\Gamma(\frac{n-1}{2})}\cdot \frac{\sqrt{\pi}}{\sqrt{c_2}}\frac{\Gamma(\frac{n-3}{2})}{\Gamma(\frac{n-2}{2})}\cdots\frac{\sqrt{\pi}}{\sqrt{c_{n-3}}}\frac{\Gamma(1)}{\Gamma(3/2)}
\int_{0}^{2\pi}\frac{d\phi_{n-2}}{c_{n-2}\cos^{2}\phi_{n-2}+c_{n-1}\sin^{2}\phi_{n-2}}
$$
Again, by \eqref{pre-integral} with $m=0$
$$
I=\frac{\sqrt{\pi}}{\sqrt{c_1}}\frac{\Gamma(\frac{n-2}{2})}{\Gamma(\frac{n-1}{2})}\cdot \frac{\sqrt{\pi}}{\sqrt{c_2}}\frac{\Gamma(\frac{n-3}{2})}{\Gamma(\frac{n-2}{2})}\cdots\frac{\sqrt{\pi}}{\sqrt{c_{n-3}}}\frac{\Gamma(1)}{\Gamma(3/2)}\cdot \frac{2\pi}{\sqrt{c_{n-2}c_{n-1}}}
=\frac{2\pi^\frac{n-1}{2}}{\Gamma(\frac{n-1}{2})}\, \prod_{i=1}^{n-1} \, c_i^{-\frac{1}{2}}.
$$
$\Box$
\vskip 4mm
\subsection{Proof of Theorem \ref{theorem:limit}}

We need several more ingredients. The first lemma is standard, see e.g., \cite{SW:1990}.   

\begin {lemma} \label{vol-diff}
Let $K$ and $L$  be a convex bodies in $\mathbb{R}^n$ such that $0 \in \text{int} (L) \subset K$. Then
$$
\vol_n(K)-\vol_n(L) = \frac{1}{n} \int_{\partial K} \langle x, N_K(x) \rangle \left[1- \left(\frac{\|x_L\|}{\|x\|}\right)^n \right] d \mu_K(x), 
$$
where $x_L=[0,x] \cap \partial L$.
\end{lemma} 
\vskip 2mm
\noindent
We can assume without loss of generality that $0 \in \text{int} (K_\delta^R)$ and when $L=K_\delta^R$, we write $x_\delta= x_{K_\delta^R}$.
\vskip 3mm
\noindent
We say that $R\, B^n_2$ is a supporting $R$-ball to $K$ if $R\, B^n_2 \cap \partial K \neq \emptyset$ and $K \subset R\, B^n_2$.
\vskip 3mm
\begin{lemma}\label{properties} 
Let $K$ be an $R$-ball convex body  and let $K_\delta^R$ be its $R$-ball floating body. 
\newline
Through every point of $\partial K^R_{\delta}$ there is at least one $R$-ball 
 that cuts off a set of volume $\delta$ exactly from $K$.
\end{lemma}
\vskip 2mm
\begin{proof}
Let $x\in \partial K^R_{\delta}$. We choose a sequence $x_{k}$, $k\in \mathbb N$,
such that $x_{k}\notin K^R_{\delta}$ and $\lim_{k\to\infty}x_{k}=x$. By the
definition of $K^R_{\delta}$ we find for every $x_{k}$ a ball $B^n_2(y_k, R)$  such
that $x_{k}\in B^n_2(y_k, R)$ and $\vol_{n} (K\setminus B^n_2(y_k, R))=\delta$. By
compactness there is a subsequence $B^n_2(y_{k_{j}}, R)$, $j \in \mathbb N$, that converges to
an $R$-ball $B^n_2(y,R)$, i.e. the sequence  of the centers converges to the center and there is a sequence of elements on $\partial B^n_2(y_{k_{j}}, R)$  that converges to
an element on  $\partial B^n_2(y, R)$.
Therefore $\text{vol}_{n}(K \setminus B^n_2(y, R))= \delta$ and $x \in  \partial B^n_2(y, R)$.
\end{proof}
\vskip 3mm
\noindent
For  the next lemma, we recall the \emph{rolling function} $r_K:\partial K \to [0,\infty)$ of a convex body $K$, which  was introduced by McMullen in \cite{MM:1974}, see also \cite{SW:1990}.
For $x \in \partial K$ with unique outer normal $N_K(x)$ it  is defined by 
\begin{equation*}
	r_K(x) = \max\{\rho: B^{n}_2 (x - \rho N_{K} (x), \rho ) \subset K\}, 
\end{equation*}
i.e., $r_K(x)$ is the maximal radius of a Euclidean ball inside $K$ that contains $x$.
If $N_K(x)$ is not unique, $r_K(x) =0$. By McMullen \cite{MM:1974} (also \cite{SW:1990}) $r_K(x)>0$  almost everywhere on $\partial K$.
It was shown in \cite{SW:1990} that for all $0 \leq \alpha <1$, 
\begin{equation}\label{r}
\int_{\partial K} \frac{1}{r_K(x)^\alpha} d \mu_K(x) < \infty.
\end{equation}
\vskip 3mm
\noindent
The next lemma and its proof is the analog of Lemma 6 of \cite{SW:1990}. 
\vskip 3mm
\begin{lemma} \label{bounded}
Let $K$ be an $R$-ball convex body in $\mathbb R^{n}$  that contains $0$ as an interior point.
Let $x \in \partial K$ such that $r_{K}(x)>0$ and let $x_{\delta}$ be the unique point in $[0,x]\cap\partial K_{\delta}^{R}$. Then there is $\delta_0$ such that for all $\delta \leq \delta_0$,  
\begin{equation}\label{bounded-1}
\frac{\|x-x_\delta\|}{\delta^\frac{2}{n+1}} \leq \gamma_n \, r_{K}(x)^{-\frac{n-1}{n+1}},
\end{equation}
where $\gamma_n$ depends on $n$ and $K$ only.
\end{lemma}
\vskip 2mm

It is important for the proof of Theorem \ref{theorem:limit} that the constant $\gamma_{n}$
only depends on the dimension $n$ and the body $K$, but not on $x\in\partial K$.
\vskip2mm

\begin{proof}
Since $0$ is an interior point of $K$ there is $s>0$ such that 
\begin{equation}\label{s-contain}
B_{2}^{n}(0,\tfrac{1}{s})\subseteq K\subseteq B_{2}^{n}(0,s).
\end{equation}
We choose $\delta_{0}=\tfrac{1}{2}\text{vol}_{n}(B_{2}^{n}(0,\tfrac{1}{4s}))$. 
In order to prove \eqref{bounded-1} we consider two cases,
 $s^{-2}r_{K}(x)\leq \|x-x_{\delta}\|_{2}$ and $s^{-2}r_{K}(x)\geq \|x-x_{\delta}\|$.
\par
\noindent
We assume first that $s^{-2}r_{K}(x)\leq \|x-x_{\delta}\|$.
More specifically, we will show 
\begin{equation}\label{FloatTh3-2}
\frac{\|x-x_\delta\|}{\delta^\frac{2}{n+1}}
\leq \frac{2^{2\frac{n-1}{n+1}} \, n^{\frac{2}{n+1}}s^{6\frac{n-1}{n+1}}}{(\text{vol}_{n-1}
(B_{2}^{n-1}))^{\frac{2}{n+1}}} \, r_{K}(x)^{-\frac{n-1}{n+1}}.
\end{equation}
Let $\delta \leq \delta_0$ and let
$z +R B^n_2$ be an $R$-ball that touches $K_\delta^R$ in $x_\delta$ and cuts off a set of volume exactly 
equal to $\delta$ from $K$,
i.e. $\operatorname{vol}_{n}(K\setminus (z +R B^n_2))=\delta$.
Such an $R$-ball exists by Theorem \ref{properties} (i). 
If $z +R B^n_2$ does not contain $B_{2}^{n}(0,\tfrac{1}{2s})$ then
 $K\setminus (z +R B^n_2)$ contains a Euclidean ball of radius
$\frac{1}{4s}$. This cannot be since $\delta$ is smaller than
$\delta_{0}=\frac{1}{2}\text{vol}_{n}(B_{2}^{n}(0,\tfrac{1}{4s}))$. Therefore
\begin{equation}\label{bounded-2}
 B_{2}^{n}(0,\tfrac{1}{2s})\subseteq z +R B^n_2.
\end{equation}
Let $T(x_\delta)$ be  the tangent hyperplane to $z +R B^n_2$ in $x_\delta$ and $T(x_\delta)^+$ the half space containing
 $z +R B^n_2$.
Then
$$
C=[x,B_{2}^{n}(0,\tfrac{1}{2s})]\cap T(x_\delta)^{-}
$$
is a cone with a base that is an ellipsoid. Instead of the cone $C$
we consider the cone 
$$
\tilde C
=[x,B_{2}^{n}(0,\tfrac{1}{2s})]\cap H(x_{\delta})^-, 
$$
where $H(x_{\delta})^{}$ is the hyperplane through $x_\delta$ and orthogonal to $x$. By a simple geometric argument
$$
\operatorname{vol}_{n}(\tilde C)\leq\operatorname{vol}_{n}(C).
$$
 The height of the  cone $\tilde C$ is $\|x-x_{\delta}\|$ and the radius of the
base is $\frac{1}{2s}\frac{\|x-x_{\delta}\|}{\|x\|}$.
Therefore 
\begin{equation*}\label{FloatTh3-1}
 \delta  = \vol_n(K \setminus (z +R\,  B^n_2)) \geq \operatorname{vol}_{n}(K\cap  H(x_{\delta})^-)
\geq \frac{1}{n}\frac{\|x-x_{\delta}\|^{n}}{\|x\|^{n-1}}\operatorname{vol}_{n-1}
(B_{2}^{n-1}(0,\tfrac{1}{2s})).
\end{equation*}
And thus
$$
\frac{\|x-x_\delta\|}{\delta^\frac{2}{n+1}} 
\leq \frac{n^{\frac{2}{n+1}} \|x\|^{2\frac{n-1}{n+1}}(2s)^{2\frac{n-1}{n+1}}}{(\text{vol}_{n-1}
(B_{2}^{n-1}))^{\frac{2}{n+1}}\|x-x_{\delta}\|^{\frac{n-1}{n+1}}}.
$$
Since $\|x\| \leq s$ by (\ref{s-contain}) and since
$s^{-2}r_{K}(x)\leq \|x-x_{\delta}\|$ by assumption
\begin{eqnarray*}
\frac{\|x-x_\delta\|}{\delta^\frac{2}{n+1}}
\leq \frac{n^{\frac{2}{n+1}} 2^{2\frac{n-1}{n+1}}s^{4\frac{n-1}{n+1}}}{(\text{vol}_{n-1}
(B_{2}^{n-1}))^{\frac{2}{n+1}}\|x-x_{\delta}\|^{\frac{n-1}{n+1}}} 
\leq \frac{n^{\frac{2}{n+1}} 2^{2\frac{n-1}{n+1}}s^{6\frac{n-1}{n+1}}}{(\text{vol}_{n-1}
(B_{2}^{n-1}))^{\frac{2}{n+1}}} r_{K}(x)^{-\frac{n-1}{n+1}}.
\end{eqnarray*}
 This proves (\ref{FloatTh3-2}).
\vskip 1mm
\noindent
Now we 
consider the case $s^{-2}r_{K}(x)\geq \|x-x_{\delta}\|$.
Again, let $T(x_\delta)$ denote the tangent hyperplane to $z +R B^n_2$ in $x_\delta$. Then, since $B_{2}^{n}(x-r_{K}(x)N_K(x),r_{K}(x))\subseteq K$
\begin{eqnarray}\label{cap}
\delta&=&
 \vol_n(K \setminus (z +R\,  B^n_2)) \geq \vol_{n}(K\cap  T(x_{\delta})^-) 
\nonumber \\
 &\geq& \vol_{n}(B_{2}^{n}(x-r_{K}(x)N_K(x),r_{K}(x))\cap T(x_{\delta})^-)  .
\end{eqnarray}
Now we are in the situation of Lemma 6 \cite{SW:1990} with $H=T(x_{\delta})$. For completeness, we give the arguments.
\par
Indeed, let $\theta$ denote the angle between $x$ and $N_{K}(x)$. 
By \eqref{s-contain} and \eqref{bounded-2} we have 
\begin{equation}\label{bound-5}
\cos\theta\geq\frac{1}{s\|x\|}\geq\frac{1}{s^{2}}.
\end{equation}
We show that
\begin{equation}\label{bound-3}
x_{\delta}\in B_{2}^{n}(x-r_{K}(x)N_K(x),r_{K}(x)).
\end{equation}
We have (the inequality follows by \eqref{bound-5})
\begin{equation}\label{bound-8}
\operatorname{vol}_{1}\left([0,x]\cap B_{2}^{n}(x-r_{K}(x)N_K(x),r_{K}(x)-{K}(x))\right)
=2 r_{K}(x)\cos\theta\geq \frac{2 r_{K}(x)}{s^{2}}.
\end{equation}
 By the assumption
$s^{-2}r_{K}(x)\geq \|x-x_{\delta}\|$ we get \eqref{bound-3}.
\par
Let $\Delta$ be the distance from $x_{\delta}$ to the boundary of
$B_{2}^{n}(x-r_{K}(x)N_K(x),r_{K}(x))$
$$
\Delta =\min\{\|x_{\delta}-z\|\ : \ \|z-(x-r_{K}(x)N_K(x))\|=r_{K}(x)\}.
$$
$B_{2}^{n}(x-r_{K}(x)N_K(x),r_{K}(x))\cap T(x_{\delta})^-$ is a cap of height $\Delta$ of the ball $B_{2}^{n}(x-r_{K}(x)N_K(x),r_{K}(x))$.
We use  e.g., Lemma 8 of  \cite{SW:1990}  to estimate the volume of this cap  
\begin{eqnarray*}
\delta &\geq&\vol_{n}(B_{2}^{n}(x-r_{K}(x)N_K(x),r_{K}(x))\cap T(x_{\delta})^-) \\
&\geq& 2(2r_{K}(x))^{\frac{n-1}{2}}\frac{\operatorname{vol}_{n-1}(B_{2}^{n-1})}{n+1}
\left\{\Delta^{\frac{n+1}{2}}-\frac{(n+1)(n-1)}{4r_{K}(x) (n+3)}
\Delta^{\frac{n+3}{2}}\right\}
\\
&=&2(2r_{K}(x))^{\frac{n-1}{2}}\frac{\operatorname{vol}_{n-1}(B_{2}^{n-1})}{n+1}
\Delta^{\frac{n+1}{2}}\left\{1-\frac{(n+1)(n-1)}{4r_{K}(x) (n+3)}
\Delta^{}\right\}.
\end{eqnarray*}
\begin{figure}[h]
			\centering
	\includegraphics[scale=0.5]{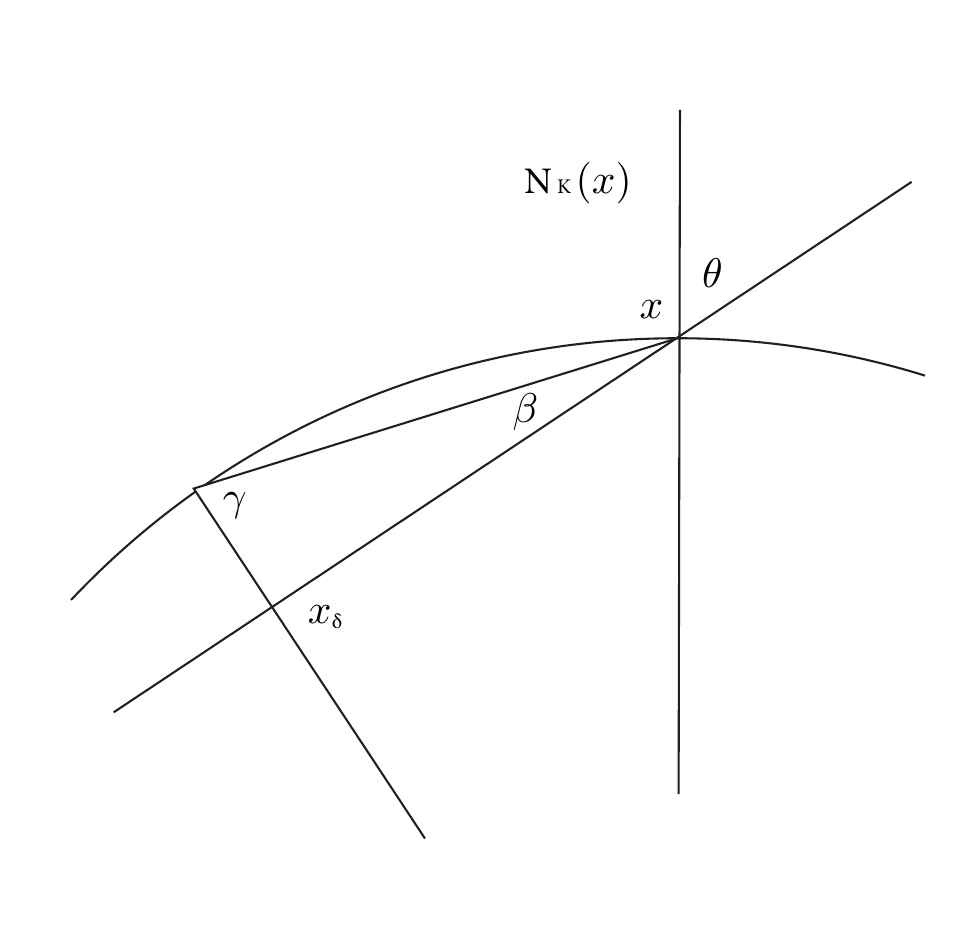}		
		\caption{The angles}
		\end{figure}	
Let $\theta$,  $\beta$  and $\gamma$ be as in Figure 1. 
If $\theta=0$ then $\Delta=\|x-x_{\delta}\|$. If $\Delta>0$, then $\gamma>0$.
By the sine law for the triangle
\begin{equation*}\label{FloatTh3-4}
\Delta=\frac{\sin \beta}{\sin\gamma}\|x-x_{\delta}\|.
\end{equation*}
Since $\beta\leq\gamma\leq\frac{\pi}{2}$ we have
$\sin\beta\leq\sin\gamma$. Therefore, $\Delta\leq\|x-x_{\delta}\|$. Moreover,
$\gamma\geq \theta$. Together with \eqref{s-contain} this implies $\sin \gamma\geq\sin\theta\geq\frac{1}{s^{2}}$. Thus
$$
\delta  \geq 
2(2r_{K}(x))^{\frac{n-1}{2}}\frac{\operatorname{vol}_{n-1}(B_{2}^{n-1})}{n+1}
(\sin\beta \|x-x_{\delta}\|)^{\frac{n+1}{2}}\left\{1-\frac{(n+1)(n-1)}{4 \, r_{K}(x) (n+3)}
  \|x-x_{\delta}\|\right\}  .
$$
It follows that
$$
r_{K}(x)^{-\frac{n-1}{n+1}}  
\geq2 \left(\frac{\operatorname{vol}_{n-1}(B_{2}^{n-1})}{n+1}\right)^{\frac{2}{n+1}}
\frac{\sin\beta \|x-x_{\delta}\|_{2}}{ \delta^{\frac{2}{n+1}}}
\left\{1-\frac{(n+1)(n-1)}{4\, r_{K}(x) (n+3)}
  \|x-x_{\delta}\|\right\}^{\frac{2}{n+1}}.
$$
Since $s^{-2}r_{K}(x)\geq\|x-x_{\delta}\|$, 
\begin{eqnarray*}
r_{K}(x)^{-\frac{n-1}{n+1}}  
&\geq& 2 \left(\frac{\vol_{n-1}(B_{2}^{n-1})}{n+1}\right)^{\frac{2}{n+1}}
\frac{\sin\beta \|x-x_{\delta}\|}{ \delta^{\frac{2}{n+1}}}
\left\{1-\frac{(n+1)(n-1)}{4\, s^{2} (n+3)}
\right\}^{\frac{2}{n+1}}, 
\end{eqnarray*}
which implies that, choosing $s \geq (n+1)^\frac{1}{2}$, 
\begin{eqnarray*}
\frac{ \|x-x_{\delta}\|}{ \delta^{\frac{2}{n+1}}} \leq \frac{1}{\sin\beta }\, \frac{1}{2} \left(\frac{n+1}{\vol_{n-1}(B_{2}^{n-1})}\right)^{\frac{2}{n+1}} \left(\frac{4}{3}\right)^{-\frac{2}{n+1}}  r(x)^{-\frac{n-1}{n+1}}.
\end{eqnarray*}
It is left to verify that $\sin\beta$ is bounded from below. 
Let $\alpha$ be the angle between $x$ and $x_{\delta}$ at $x-r_{K}(x)N_K(x)$.
We have
\begin{equation}\label{bound-9}
\gamma=\beta+\theta
\hskip 20mm
2\gamma+\alpha=\pi
\hskip 20mm
\alpha+\theta+\frac{\pi}{2}\leq\pi
\end{equation}
Indeed, $\gamma=\beta+\theta$ follows immediately from Figure 1. For $2\gamma+\alpha=\pi$
we use that the angle sum in a triangle equals $\pi$. By \eqref{bound-8} and by the assumption
$s^{-2}r_{K}(x)\geq \|x-x_{\delta}\|$ we conclude that
the angle at $x_{\delta}$ between $x$ and $x-r_{K}(x)N_K(x)$ is greater than $\frac{\pi}{2}$.
This implies $\alpha+\theta+\frac{\pi}{2}\leq\pi$.
\par
By \eqref{bound-9}
$$
\tfrac{1}{2}(\tfrac{\pi}{2}-\theta)\leq\beta.
$$
On the interval $[0,\frac{\pi}{2}]$ the function $\sin$
is increasing. By \eqref{bound-5}
$$
\sin(2\beta)\geq\sin(\tfrac{\pi}{2}-\theta)
=\sin(\tfrac{\pi}{2})\cos(-\theta)
+\cos(\tfrac{\pi}{2})\sin(-\theta)
=\cos\theta\geq\tfrac{1}{s^{2}}.
$$
\end{proof}
\vskip 3mm
\noindent
In the next lemma we estimate the volume of a cap $C(h, \mathcal{E})$ of height $h$ cut off  from an ellipsoid by an $R$-ball.
More precisely, let an ellipsoid $\mathcal{E}$  be given by 
\begin{equation}\label{ellipse}
\sum_{i=1}^{n-1} \frac{x_i^2}{a_i^2} + \frac{(x_n-a_n)^2}{a_n^2}=1
\end{equation}
and let the $R$-ball have equation
\begin{equation}\label{cut-ball}
\sum_{i=1}^{n-1}x_i^2 + (x_n-a)^2=R^2,
\end{equation}
i.e., the $R$-ball is centered at $z=(0 \cdots, 0,a)$.
Without loss of generality, we can assume that $a > R \geq a_n \geq a_{n-1} \dots \geq a_1$.
Then  
$$
C(h, \mathcal{E}) = \mathcal{E} \setminus B^n_2(z,R)
$$
and the height $h$ of the cap is $h=a-R$. 
\begin {lemma} \label{cap} Let $\varepsilon >0$ be given. Let the ellipsoid $\mathcal{E}$ be such that $\frac{a_n}{{a_i}^2} > \frac{1}{R}$ for  at least one $i$ . Then we have for sufficiently small $h$, 
\begin{eqnarray*}
&\hskip -10mm (1 - d \varepsilon) \frac{2^\frac{n+1}{2}}{(n-1)(n+1)} (a-R)^\frac{n+1}{2}\,  \int_{S^{n-2}} \frac{d\sigma(\xi)}
{\left(\sum_{i=1}^{n-1} \left(\frac{a_n}{a_i^2} -\frac{1}{R} \right) \xi_i^2\right)^\frac{n-1}{2}}
\leq 
\vol_n(C(h, \mathcal{E}))  \\
&\hskip 10mm \leq (1 + c \varepsilon) \frac{2^\frac{n+1}{2}}{(n-1)(n+1)} (a-R)^\frac{n+1}{2}\,  \int_{S^{n-2}} \frac{d\sigma(\xi)}
{\left(\sum_{i=1}^{n-1} \left(\frac{a_n}{a_i^2} -\frac{1}{R} \right) \xi_i^2\right)^\frac{n-1}{2}},
\end{eqnarray*} 
where $c$ and $d$ are constants.
\end{lemma}
\vskip 2mm
\begin{proof}
We put 
$$
f(x) = a_n -a_n \left(1- \sum_{i=1}^{n-1} \frac{x_i^2}{a_i^2}  \right)^\frac{1}{2}
$$
and
$$
g(x) = a - R \left(1- \sum_{i=1}^{n-1} \frac{x_i^2}{R^2}  \right)^\frac{1}{2}.
$$
Passing to polar coordinates, those become
$$
f(\rho, \xi) = a_n -a_n \left(1- \rho^2 \sum_{i=1}^{n-1} \frac{\xi_i^2}{a_i^2}  \right)^\frac{1}{2},  \hskip 4mm g(\rho) = a - R \left(1-  \frac{\rho^2}{R^2}  \right)^\frac{1}{2}.
$$
The intersection of the two functions has equation
$$
a-a_n = R \left(1-  \frac{\rho^2}{R^2}  \right)^\frac{1}{2} - a_n \left(1- \rho^2 \sum_{i=1}^{n-1} \frac{\xi_i^2}{a_i^2}  \right)^\frac{1}{2}.
$$
We will apply the lemma for $\rho$ small and therefore estimate with  suitable constants $c$ and $d$,
$$
a-a_n \leq  R \left(1-  \frac{1}{2} \,  \frac{\rho^2}{R^2}  \right) - a_n \left(1- \frac{\rho^2}{2} \sum_{i=1}^{n-1} \frac{\xi_i^2}{a_i^2} (1+c \rho^2) \right),
$$
$$
a-a_n \geq  R \left(1-  \frac{1}{2} \,  \frac{\rho^2}{R^2} (1+d \rho^2)  \right) - a_n \left(1- \frac{\rho^2}{2} \sum_{i=1}^{n-1} \frac{\xi_i^2}{a_i^2}  \right)
$$
or, equivalently, with new constants $c$ and $d$,
$$
 \frac{\rho^2}{2} \, \left(a_n  \sum_{i=1}^{n-1} \frac{\xi_i^2}{a_i^2}  -\frac{1}{R} \right) +d   \rho^4 \leq a- R \leq  \frac{\rho^2}{2} \, \left(a_n  \sum_{i=1}^{n-1} \frac{\xi_i^2}{a_i^2}  -\frac{1}{R} \right) +c \rho^4.
$$
This leads to, for $\varepsilon >0$ given, for sufficiently small $\rho$, 
\begin{equation}\label{rho}
\sqrt{2} (1-d \varepsilon) \frac{(a-R)^\frac{1}{2}}{\left(a_n  \sum_{i=1}^{n-1} \frac{\xi_i^2}{a_i^2}  -\frac{1}{R} \right)^\frac{1}{2}} \leq \rho= \rho(\xi) \leq \sqrt{2} (1+c \varepsilon) \frac{(a-R)^\frac{1}{2}}{\left(a_n  \sum_{i=1}^{n-1} \frac{\xi_i^2}{a_i^2}  -\frac{1}{R} \right)^\frac{1}{2}}.
\end{equation}
Then, using again that $(1-x^2)^\frac{1}{2} \leq 1 - \frac{x^2}{2}$, 
\begin{eqnarray*}
\vol_n(C(h, \mathcal{E})) &=& \int_{S^{n-2}}  \int_{\rho =0}^{\rho(\xi)} \left(g(\rho) - f(\rho, \xi)\right) \rho^{n-2} d\rho\,  d\sigma(\xi) \\
&\leq& 
(1 + c \varepsilon) \frac{2^\frac{n+1}{2}}{(n-1)(n+1)} (a-R)^\frac{n+1}{2}\,  \int_{S^{n-2}} \frac{d\sigma(\xi)}
{\left(\sum_{i=1}^{n-1} \left(\frac{a_n}{a_i^2} -\frac{1}{R} \right) \xi_i^2\right)^\frac{n-1}{2}}
\end{eqnarray*}
and similarly for the lower bound.
\end{proof}
\vskip 3mm

\begin {lemma} \label{limit} 
Let $K$ be an $R$-ball convex body in $\mathbb R^{n}$.
Let $x \in \partial K$ be such that $\kappa_i(x, K) > \frac{1
}{R}$ for all $i=1,\dots,n-1$. Then
$$
\lim_{\delta \to 0} \left\langle \frac{x}{\|x\|}, N_K(x) \right\rangle  \frac{ \|x-x_\delta\|}{\delta^\frac{2}{n+1}} =  
\frac{1}{2} \left(n^{2}-1\right)^\frac{2}{n+1} 
\left[\int_{S^{n-2}} \frac{d\sigma(\xi)}
{\left(\sum_{i=1}^{n-1} \left(\kappa_i(x, K) -\frac{1}{R} \right) \xi_i^2\right)^\frac{n-1}{2}}\right]^{-\frac{2}{n+1}}.
$$
\end{lemma}
\vskip 4mm
\noindent
For the proof of Lemma \ref{limit}, we need more ingredients (see e.g., \cite{SW4}).
\par
\noindent
Let $K$ be a convex body in $\mathbb R^{n}$ such that $0\in\partial K$
and $N_{K}(0)=e_{n}$. Moreover suppose that $\sum_{i=1}^{n-1}\frac{\xi_{i}^{2}}{b_{i}^{2}}=1$
is the indicatrix of Dupin. Then the principal Gauss-Kronecker curvatures are $b_{i}^{2}$,
$i=1,\dots,n-1$ and the Gauss-Kronecker curvature of $K$ at $0$ is $\prod_{i=1}^{n-1}b_{i}^{2}$.
We call $\mathcal E$ the standard approximating ellipsoid, where
\begin{equation}\label{StandardEll1}
\mathcal E
=\left\{\xi\in\mathbb R^{n}\left|
\sum_{i=1}^{n-1}\frac{\xi_{i}^{2}}{b_{i}^{2}}+
\frac{\left(\xi_{n}-\left(\prod_{i=1}^{n-1}b_{i}\right)^{\frac{2}{n-1}}\right)^{2}}
{\left(\prod_{i=1}^{n-1}b_{i}\right)^{\frac{2}{n-1}}}
\leq\left(\prod_{i=1}^{n-1}b_{i}\right)^{\frac{2}{n-1}}
\right.\right\}.
\end{equation}
There is an increasing, continuous function 
$\phi:[0,\infty)\to[1,\infty)$ with
$\phi(0)=1$ such that for all sufficiently small $t\geq0$
\begin{eqnarray}\label{ApproxEllip1-11}
&&\left\{\left.\left(\frac{\xi_{1}}{\phi(t)},\dots,\frac{\xi_{n-1}}{\phi(t)},t\right)\right|
\xi\in\mathcal E\right\}
\\
&&\subseteq K\cap H((0,\dots,0,t),N(0))
\subseteq
\left\{(\phi(t)\xi_{1},\dots,\phi(t)\xi_{n-1},t)|\xi\in\mathcal E\right\}.
\nonumber
\end{eqnarray}
We put $b_{i}=\frac{a_{i}}{\sqrt{a_{n}}}$ and $a_{n}=\left(\prod_{i=1}^{n-1}b_{i}\right)^{\frac{2}{n-1}}$
and we get for \eqref{StandardEll1}
\begin{equation}\label{StandardEll2}
\mathcal E
=\left\{\xi\in\mathbb R^{n}\left|
\sum_{i=1}^{n-1}\frac{a_{n}}{a_{i}^{2}}\xi_{i}^{2}+
\frac{\left(\xi_{n}-a_{n}\right)^{2}}
{a_{n}}
\leq a_{n}
\right.\right\}.
\end{equation}
This means
$$
\xi_{n}=a_{n}-a_{n}\sqrt{1-\sum_{i=1}^{n-1}\frac{1}{a_{i}^{2}}\xi_{i}^{2}}.
$$
We will need  two further ellipsoids $\mathcal{E} (\varepsilon^-)$
and $\mathcal{E} (\varepsilon^+)$, one  slightly smaller than $\mathcal E$
and the other slightly bigger.
\begin{equation}\label{StandardEll4}
\mathcal E(\varepsilon^-)
=\left\{\xi\in\mathbb R^{n}\left|
\sum_{i=1}^{n-1}\frac{a_{n}}{(1-\varepsilon)^{2}a_{i}^{2}}\xi_{i}^{2}+
\frac{\left(\xi_{n}-a_{n}\right)^{2}}
{a_{n}}
\leq a_{n}
\right.\right\}.
\end{equation}
\begin{equation}\label{StandardEll5}
\mathcal E(\varepsilon^+)
=\left\{\xi\in\mathbb R^{n}\left|
\sum_{i=1}^{n-1}\frac{a_{n}}{(1+\varepsilon)^{2}a_{i}^{2}}\xi_{i}^{2}+
\frac{\left(\xi_{n}-a_{n}\right)^{2}}
{a_{n}}
\leq a_{n}
\right.\right\}.
\end{equation}
\vskip3mm
\noindent
Now we give the proof of Lemma \ref{limit}.
\par
\noindent
\begin{proof}
Let $0$ be an interior point of $K$ and let
$x \in  \partial K$ with outer normal $N_K(x)$. We can assume $N_K(x)=e_n$. Then
locally around $x$, $\partial K$  can be approximated by the ellipsoid $x+\mathcal E$,
where $\mathcal E$ is given by \eqref{StandardEll2}. 
\newline
We can assume that the axes of principal curvature are the unit vectors $e_{1},\dots,e_{n-1}$.
Let $\varepsilon >0$  be given.
We can assume that 
$$
R \geq a_n \geq \dots \geq a_1.
$$
Let $\mathcal{E} (\varepsilon^-)$ be the ellipsoid centered at $ a_{n} e_n$ 
whose principal axes  coincide with the ones of $\mathcal{E}$,  but have lengths $(1-\varepsilon) a_1, \dots, (1-\varepsilon) a_{n-1}, a_{n}$. Similarly, let $\mathcal{E} (\varepsilon^+)$ be the ellipsoid centered at $a_{n} e_n$, with the same principal axes as $\mathcal{E}$,  but with lengths $(1+\varepsilon) a_1, \dots, (1+\varepsilon) a_{n-1}, a_{n}$, i.e., 
the ellipsoids given in (\ref{StandardEll4}) resp. (\ref{StandardEll5}).
Then
$$ 
x \in \partial \mathcal{E}  \hskip 5mm \text { and } \hskip 5mm N_{\mathcal{E}}(x) = N_K(x),
$$
and (see, e.g., \cite{SW4})  there exists a $\Delta_\varepsilon >0$ such that 
\begin{eqnarray}\label{ellipse}
 H^-\left(x- \Delta_\varepsilon e_n, e_n\right)  \  \cap  \  \mathcal{E} (\varepsilon^-)
  \subseteq  H^-\left(x- \Delta_\varepsilon e_n, e_n\right) \  \cap \   K 
\subseteq H^-\left(x- \Delta_\varepsilon e_n, e_n\right) \  \cap  \  \mathcal{E} (\varepsilon^+) .
\end{eqnarray}
Let $x_\delta \in \partial K_\delta^R\cap[0,x]$. We choose $\delta$ so small that  for all support balls $z +R B^n_2$ to $K^R_\delta$  in $x_\delta$ we have
\begin{equation}\label{delta}
 \mathcal{E} (\varepsilon^-)\setminus  (z +R B^n_2) \subseteq   H^-\left( x-\Delta_\varepsilon e_n, e_n\right) \,  \cap  \,  \mathcal{E} (\varepsilon^-) .
 \end{equation}
Let
$H\left(x_\delta, e_n\right)$ be the hyperplane through $x_\delta$ and orthogonal to $e_n$. Then
$$
x-\left\langle \frac{x}{\|x\|}, N_K (x)\right\rangle  \|x-x_\delta\|  e_n 
$$ 
is the intersection of $H\left(x_\delta, e_n\right)$ with the line
$x+\{\lambda e_{n}|\lambda\in\mathbb R\}$.
Since $x_\delta \in \partial K_\delta^R$ and since the  $R$-ball centered at $a=x-\left(R+ \left\langle \frac{x}{\|x\|}, N_K (x)\right\rangle \|x-x_\delta\| \right) e_n$ with radius $R$ does not contain $x_{\delta}$ 
$$
\delta \leq \vol_n\left(\mathcal{E} (\varepsilon^+) \setminus B^n_2\left(a, R\right) \right).
$$
Thus, with Lemma \ref{cap},
$$
\delta \leq 
 \frac{(1 + c \varepsilon) \, 2^\frac{n+1}{2}}{(n-1)(n+1)} \left[\left\langle \frac{x}{\|x\|}, N_K (x)\right\rangle \|x-x_\delta\|\right]^\frac{n+1}{2}\,  \int_{S^{n-2}} \frac{d\sigma(\xi)}
{\left(\sum_{i=1}^{n-1} \left(\frac{a_n}{(1+\varepsilon)^2 a_i^2} -\frac{1}{R} \right) \xi_i^2\right)^\frac{n-1}{2}}.
$$
Now we treat the estimate from below. We keep the same coordinate setup as above. In particular,   $x\in\partial K$ with outer normal $N_K(x)=e_n$.
Let $\delta >0$ be so small that (\ref{delta}) holds and let $x_\delta \in \partial K_\delta$.
We denote  by  $\theta$ the angle between
$e_n$  and $x$. We show
\begin{eqnarray}\label{Limit1}
&&\hskip -27mm \left|1- 
\frac{a_n^{3}\sin\theta^{2}\|x-x_{\delta}\|^{2}(1 - \varepsilon)a_1)^2}{\left((1 - \varepsilon)a_1)^2-a_{n}
\cos^{2}\theta\|x-x_{\delta}\| \right)^{2}}\right|    \cos\theta \|x_{\delta}-x\|\leq  \nonumber\\
&&\hskip 30mm\leq\inf_{y\in\partial \mathcal{E}(\varepsilon^-)}\|x_{\delta}-y\|
\leq\cos\theta \|x_{\delta}-x\|.
\end{eqnarray}
The right hand side inequality is obvious. We show the left hand side inequality.
The left hand side inequality reduces immediately to a $2$-dimensional inequality.
The expression $\inf_{y\in\partial \mathcal{E}(\varepsilon^-)}\|x_{\delta}-y\|$ is smallest if the vector $\|x- x_\delta\| (sin \theta, 0)$ is parallel to the direction 
of the principal axis of $\mathcal{E}$ that has length $a_1$.
\par
After shifting $x$ to $0$, the vector $x_\delta$ can be written
in the $2$-dimensional coordinate system spanned by $e_n$ and $e_{1}$, 
as
$$
x_\delta = \|x- x_\delta\| (\text{sin}  \theta, \text{cos} \theta).
$$ 
The $2$-dimensional ellipsoid is given by
$$
\xi(2)=a_{n}-a_{n}\sqrt{1-\frac{\xi(1)^{2}}{((1-\varepsilon)a_{1})^{2}}}.
$$
Now we compute the distance of $x_{\delta}$ to $\mathcal{E}(\varepsilon^-)$. We find the minimum of
$$
\|(x_{\delta}(1),x_{\delta}(2))-(\xi(1),\xi(2))\|^{2}
=|x_{\delta}(1)-\xi(1)|^{2}+\left|x_{\delta}(2)-a_{n}-a_{n}\sqrt{1-\frac{\xi(1)^{2}}{a((1-\varepsilon)a_{1})^{2}}} \right|^{2}
$$
by differentiating with respect to $\xi(1)$.
The derivative of this expression is
\begin{eqnarray*}
&&2(\xi(1)-x_{\delta}(1))
+2\left(x_{\delta}(2)-a_{n}-a_{n}\sqrt{1-\frac{\xi(1)^{2}}{((1-\varepsilon)a_{1})^{2}}} \right)\frac{a_{n}\xi(1)}{((1-\varepsilon)a_{1})^{2}\sqrt{1-\frac{\xi(1)^{2}}{((1-\varepsilon)a_{1})^{2}}}}
\\
&&=2(\xi(1)-x_{\delta}(1))
+2\frac{a_{n}\xi(1)}{a_{1}^{2}}\left(\frac{x_{\delta}(2)-a_{n}}{\sqrt{1-\frac{\xi(1)^{2}}{((1-\varepsilon)a_{1})^{2}}} }-a_{n}\right)
\end{eqnarray*}
For the minium we get
\begin{eqnarray*}
x_{\delta}(1)
&=&\xi(1)\left(1-
\frac{a_{n}}{((1-\varepsilon)a_{1})^{2}}\left(\frac{x_{\delta}(2)-a_{n}}{\sqrt{1-\frac{\xi(1)^{2}}{((1-\varepsilon)a_{1})^{2}}} }-a_{n}\right)
\right)
\\
&\geq&\xi(1)\left(1-
\frac{a_{n}}{((1-\varepsilon)a_{1})^{2}}\left(\frac{x_{\delta}(2)-a_{n}}{1-\frac{\xi(1)^{2}}{((1-\varepsilon)a_{1})^{2}} }-a_{n}\right)
\right)
\\
&=&\xi(1)\left(1-
\frac{a_{n}}{((1-\varepsilon)a_{1})^{2}}\left(x_{\delta}(2)+\frac{(x_{\delta}(2)-a_{n})\frac{\xi(1)^{2}}{((1-\varepsilon)a_{1})^{2}}}{1-\frac{\xi(1)^{2}}{((1-\varepsilon)a_{1})^{2}} }\right)
\right)
\end{eqnarray*}
We may assume that $\xi(1)^{2}\leq 2((1-\varepsilon)a_{1})^{2}$. Then
\begin{eqnarray*}
x_{\delta}(1)
&\geq&\xi(1)\left(1-
\frac{a_{n}}{((1-\varepsilon)a_{1})}^{2}x_{\delta}(2)\right)+2(x_{\delta}(2)-a_{n})\frac{a_{n}\xi(1)^{3}}{((1-\varepsilon)a_{1})^{4}}
\end{eqnarray*}
We may also assume that $a_{n}\geq x_{\delta}(2)$. Therefore
$$
x_{\delta}(1)\geq \xi(1)\left(1-
\frac{a_{n}}{((1-\varepsilon)a_{1})^{2}}x_{\delta}(2)\right).
$$
Therefore
\begin{eqnarray*}
&&\|(x_{\delta}(1),x_{\delta}(2))-(\xi(1),\xi(2)\|^{2}
=\left\|(x_{\delta}(1),x_{\delta}(2))-\left(\xi(1),a_{n}-a_{n}\sqrt{1-\frac{\xi(1)^{2}}{((1-\varepsilon)a_{1})^{2}}}\right)\right\|^{2}
\\
&&\geq
\left|x_{\delta}(2)-a_{n}+a_{n}\sqrt{1-\frac{\xi(1)^{2}}{((1-\varepsilon)a_{1})^{2}}}\right|^{2}
\geq\left|x_{\delta}(2)-a_{n}+a_{n}\left(1-\frac{\xi(1)^{2}}{((1-\varepsilon)a_{1})^{2}}\right)\right|^{2}
\\
&&
\\
&&=\left|x_{\delta}(2)-\frac{a_{n}\xi(1)^{2}}{((1-\varepsilon)a_{1})^{2}}\right|^{2}
\geq\left|x_{\delta}(2)-\frac{a_{n}x_{\delta}(1)^{2}}{((1-\varepsilon)a_{1})^{2}\left(1-
\frac{a_{n}}{((1-\varepsilon)a_{1})^{2}}x_{\delta}(2)\right)^{2}}\right|^{2}.
\end{eqnarray*}
The  $R$-ball centered at $a=( R+ d_0) e_n$ with 
$d_{0}=\inf_{y\in\partial \mathcal{E}(\varepsilon^-)}\|x_{\delta}-y\|$ and with radius $R$ through $d_0 e_n$ cuts off strictly less 
than $\delta$ from $K$  as $x_\delta$ is in the interior of this $R$-ball.
Thus, with Lemma \ref{cap}, and observing that $\text{cos} \, \theta = \langle \frac{x}{\|x\|}, N_K (x)\rangle$, we get with (new)  constants $c$ and $d$
\begin{eqnarray*}
\delta &\geq& \vol_n\left(\mathcal{E} (\varepsilon^-) \setminus B^n_2\left(a, R\right) \right) \\
&\geq &
 \frac{(1 -d \varepsilon) \, 2^\frac{n+1}{2}}{(n-1)(n+1)} (1-c \varepsilon)^\frac{n+1}{2} \left[\left\langle \frac{x}{\|x\|}, N_K (x)\right\rangle \|x-x_\delta\|\right]^\frac{n+1}{2}\,  
 \\
 &&\hskip20mm
 \int_{S^{n-2}} \frac{d\sigma(\xi)}
{\left(\sum_{i=1}^{n-1} \left(\frac{a_n}{(1-\varepsilon)^2 a_i^2} -\frac{1}{R} \right) \xi_i^2\right)^\frac{n-1}{2}}.
\end{eqnarray*}
\end{proof}
\vskip 4mm
\noindent
Now we are ready to  prove Theorem~\ref{theorem:limit}.
\vskip 2mm
\noindent
{\bf Proof of Theorem \ref{theorem:limit}}
\vskip 2mm
\noindent
We put $x_\delta=[0,x] \cap \partial K_\delta^R$. Then we have by Lemma \ref{vol-diff},
$$
\vol_n(K)-\vol_n( K_\delta^R) = \frac{1}{n} \int_{\partial K} \langle x, N_K(x) \rangle \left[1- \left(\frac{\|x_\delta\|}{\|x\|}\right)^n \right] d \mu_K(x). 
$$
For $\delta$ small enough, we estimate
\begin{align*}
  &n \, \langle \frac{x}{\|x\|}, N_K(x)  \rangle \|x-x_\delta\| \, (1 - c \|x-x_\delta\|)\\
&\leq \langle x, N_K(x) \rangle \left[1- \left(\frac{\|x_\delta\|}{\|x\|}\right)^n \right] = \langle x, N_K(x) \rangle \left[1- \left(1-\frac{\|x-x_\delta\|}{\|x\|}\right)^n \right]\\
&\leq n \, \left\langle \frac{x}{\|x\|}, N_K(x)  \right\rangle \|x-x_\delta\|,
\end{align*}
where $c$ is an absolute constant.
Thus
\begin{align*}
&\lim_{\delta \to 0} \int_{\partial K} \left\langle \frac{x}{\|x\|}, N_K(x) \right\rangle \frac{ \|x-x_\delta\| (1 - c \|x-x_\delta\|)}{\delta^\frac{2}{n+1}} d \mu_K(x)\\
&\leq \lim_{\delta \to 0} \frac{\vol_n(K) - \vol_n(K_\delta^R)} {\delta^\frac{2}{n+1}}  \\
&\leq \lim_{\delta \to 0} \int_{\partial K} \left\langle \frac{x}{\|x\|}, N_K(x) \right\rangle  \frac{ \|x-x_\delta\|}{\delta^\frac{2}{n+1}} d \mu_K(x).
\end{align*}
By (\ref{r}) and Lemma \ref{bounded} and  we can interchange integration and limit. Then with Lemma \ref{limit} and Proposition \ref{integral},  we obtain that
\begin{align*}
\lim_{\delta \to 0} \frac{\vol_n(K) - \vol_n(K_\delta^R)} {\delta^\frac{2}{n+1}}  = c_n 
\int_{\partial K}  \prod _{i=1}^{n-1} \left( \kappa_i (K, x) -\frac{1}{R} \right)^\frac{1}{n+1} d \mu_K(x).
\end{align*}

%\end{proof}

\subsection{Proof of Proposition \ref{rel-asa-properties}}
%\begin{proof} 
\vskip 2mm
\noindent
{\bf Proof of (i)}
\par
\noindent
$K$ is $R$-ball convex is equivalent to $a K$ is $a R$-ball convex. Thus the proof follows immediately from the definition or  from Theorem \ref{theorem:limit} and (\ref{Affine:map:floating:body}).

\vskip 3mm
\noindent
{\bf Proof of (ii)}

	\noindent Let  $C$ and $K$ be $R$-convex bodies such that $C\cup K$ is $R$-convex. Note that $C \cap K$ is $R$-convex.
	To prove the valuation property, we follow the approach of \cite{Sc}.
	We decompose
	\begin{eqnarray*}
		\partial(C\cup K)
		&=&\{\partial C\cap\partial K\}\cup\{\partial C\cap K^{c}\}\cup\{C^{c}\cap \partial K\}   \\
		\partial(C\cap K)
		&=&\{\partial C\cap\partial K\}
		\cup\{\partial C\cap\operatorname{int}(K)\}\cup\{\operatorname{int}(C)\cap \partial K\}    \\
		\partial C&=&\{\partial C\cap\partial K\}\cup\{\partial C\cap K^{c}\}
		\cup\{\partial C\cap\operatorname{int}(K)\}   \\
		\partial K&=&\{\partial C\cap\partial K\}\cup\{\partial K\cap C^{c}\}
		\cup\{\partial K\cap\operatorname{int}(C)\}.
	\end{eqnarray*}
We have to show
\begin{eqnarray*}
&& \int_{\partial (K \cup C)} \prod _{i=1}^{n-1} \left( \kappa_i (K\cup C, x) -\frac{1}{R} \right)^\frac{1}{n+1} d \mu_{K\cup C}(x) 
+\\	  &&\int_{\partial (K \cap C) } \prod _{i=1}^{n-1} \left( \kappa_i (K\cap C, x) -\frac{1}{R} \right)^\frac{1}{n+1} d \mu_{K\cap C}(x)\\
 &&=  \int_{\partial K } \prod _{i=1}^{n-1} \left( \kappa_i (K, x) -\frac{1}{R} \right)^\frac{1}{n+1} d \mu_{K}(x) +  \int_{\partial C } \prod _{i=1}^{n-1} \left( \kappa_i (C, x) -\frac{1}{R} \right)^\frac{1}{n+1} d \mu_{ C}(x).
\end{eqnarray*}	
By the above decomposition
\begin{eqnarray}\label{Valu1}
&& \int_{\partial (K \cup C)} \prod _{i=1}^{n-1} \left( \kappa_i (K\cup C, x) -\frac{1}{R} \right)^\frac{1}{n+1} d \mu_{K\cup C}(x)	 =  \nonumber\\
&& \int_{\partial K \cap \partial C} \prod _{i=1}^{n-1} \left( \kappa_i (K\cup C, x) -\frac{1}{R} \right)^\frac{1}{n+1} d \mu_{K\cup C}(x)	 + \int_{\partial C \cap K^c } \prod _{i=1}^{n-1} \left( \kappa_i (K\cup C, x) -\frac{1}{R} \right)^\frac{1}{n+1} d \mu_{K\cup C}(x)+ \nonumber \\
&&\int_{\partial K \cap \partial C^c } \prod _{i=1}^{n-1} \left( \kappa_i (K\cup C, x) -\frac{1}{R} \right)^\frac{1}{n+1} d \mu_{K\cup C}(x)	=  \nonumber\\
&&  \int_{\partial K \cap \partial C} \prod _{i=1}^{n-1} \left( \kappa_i (K\cup C, x) -\frac{1}{R} \right)^\frac{1}{n+1} d \mu_{K\cup C}(x)	 + \int_{\partial C \cap K^c } \prod _{i=1}^{n-1} \left( \kappa_i (C, x) -\frac{1}{R} \right)^\frac{1}{n+1} d \mu_{C}(x)+  \nonumber\\
&&\int_{\partial K \cap C^c } \prod _{i=1}^{n-1} \left( \kappa_i (K, x) -\frac{1}{R} \right)^\frac{1}{n+1} d \mu_{K}(x)	
\end{eqnarray}
and
\begin{eqnarray}\label{Valu2}
&& \int_{\partial (K \cap C)} \prod _{i=1}^{n-1} \left( \kappa_i (K\cap C, x) -\frac{1}{R} \right)^\frac{1}{n+1} d \mu_{K\cap C}(x)	 = \nonumber\\
&& \int_{\partial K \cap \partial C} \prod _{i=1}^{n-1} \left( \kappa_i (K\cap C, x) -\frac{1}{R} \right)^\frac{1}{n+1} d \mu_{K\cap C}(x)	 + \int_{\partial C \cap \text{int}K } \prod _{i=1}^{n-1} \left( \kappa_i (K\cap C, x) -\frac{1}{R} \right)^\frac{1}{n+1} d \mu_{K\cap C}(x)+ \nonumber\\
&&\int_{\partial K \cap  \text{int}C } \prod _{i=1}^{n-1} \left( \kappa_i (K\cap C, x) -\frac{1}{R} \right)^\frac{1}{n+1} d \mu_{K\cap C}(x)	= \nonumber\\
&&  \int_{\partial K \cap \partial C} \prod _{i=1}^{n-1} \left( \kappa_i (K\cap C, x) -\frac{1}{R} \right)^\frac{1}{n+1} d \mu_{K\cap C}(x)	 + \int_{\partial C \cap \text{int} K } \prod _{i=1}^{n-1} \left( \kappa_i (C, x) -\frac{1}{R} \right)^\frac{1}{n+1} d \mu_{C}(x)+ \nonumber\\
&&\int_{\partial K \cap \text{int}C } \prod _{i=1}^{n-1} \left( \kappa_i (K, x) -\frac{1}{R} \right)^\frac{1}{n+1} d \mu_{K}(x).	
\end{eqnarray}
Adding (\ref{Valu1}) and (\ref{Valu2}), it remains to show that 		
\begin{eqnarray*}
&&\int_{\partial K \cap \partial C} \prod _{i=1}^{n-1} \left( \kappa_i (K\cup C, x) -\frac{1}{R} \right)^\frac{1}{n+1} d \mu_{K\cup C}(x) + \int_{\partial K \cap \partial C} \prod _{i=1}^{n-1} \left( \kappa_i (K\cap C, x) -\frac{1}{R} \right)^\frac{1}{n+1} d \mu_{K\cap C}(x)	\\
&&= \int_{\partial K \cap \partial C } \prod _{i=1}^{n-1} \left( \kappa_i (K, x) -\frac{1}{R} \right)^\frac{1}{n+1} d \mu_{K}(x)	+ 
\int_{\partial K \cap \partial C } \prod _{i=1}^{n-1} \left( \kappa_i (C, x) -\frac{1}{R} \right)^\frac{1}{n+1} d \mu_{C}(x).
\end{eqnarray*}
Please note that $d\mu_{C\cup K}=d\mu_{C\cap K}$ on $\partial K\cap \partial C$. 
	This holds because both measure are equal to the $(n-1)$-dimensional
	Hausdorff measure.  Therefore
\begin{eqnarray*}
&&\int_{\partial K \cap \partial C} \prod _{i=1}^{n-1} \left( \kappa_i (K\cup C, x) -\frac{1}{R} \right)^\frac{1}{n+1} d \mu_{K\cup C}(x) + \int_{\partial K \cap \partial C} \prod _{i=1}^{n-1} \left( \kappa_i (K\cap C, x) -\frac{1}{R} \right)^\frac{1}{n+1} d \mu_{K\cap C}(x)	\\
&&=\int_{\partial K \cap \partial C} \prod _{i=1}^{n-1} \left( \kappa_i (K\cup C, x) -\frac{1}{R} \right)^\frac{1}{n+1} d \mu_{K\cap C}(x) + \int_{\partial K \cap \partial C} \prod _{i=1}^{n-1} \left( \kappa_i (K\cap C, x) -\frac{1}{R} \right)^\frac{1}{n+1} d \mu_{K\cap C}(x)	\\
&&=\int_{\partial K \cap \partial C} \prod _{i=1}^{n-1}\left( \text{min}\left\{\kappa_i (K, x), \kappa_i(C,x)\right\} -\frac{1}{R} \right)^\frac{1}{n+1} d \mu_{K\cap C}(x) \\
&&+ \int_{\partial K \cap \partial C} \prod _{i=1}^{n-1} \left(\text{max}\left\{\kappa_i (K, x), \kappa_i(C,x)\right\}  -\frac{1}{R} \right)^\frac{1}{n+1} d \mu_{K\cap C}(x)	\\
&&=\int_{\{x \in \partial K \cap \partial C:  \text{min}\left\{\kappa_i (K, x), \kappa_i(C,x)\right\}=  \kappa_i(C,x)\}} \prod _{i=1}^{n-1}\left(\kappa_i(C,x)-\frac{1}{R} \right)^\frac{1}{n+1} d \mu_{K\cap C}(x) +  \\
&&+ \int_{\{x \in \partial K \cap \partial C:  \text{min}\left\{\kappa_i (K, x), \kappa_i(C,x)\right\}=  \kappa_i(K,x)\}} \prod _{i=1}^{n-1}\left( \kappa_i (K, x) -\frac{1}{R} \right)^\frac{1}{n+1} d \mu_{K\cap C}(x)\\
&&+\int_{\{x \in \partial K \cap \partial C:  \text{max}\left\{\kappa_i (K, x), \kappa_i(C,x)\right\}=  \kappa_i(C,x)\}} \prod _{i=1}^{n-1}\left(\kappa_i(C,x)-\frac{1}{R} \right)^\frac{1}{n+1} d \mu_{K\cap C}(x) +  \\
&&+ \int_{\{x \in \partial K \cap \partial C:  \text{max}\left\{\kappa_i (K, x), \kappa_i(C,x)\right\}=  \kappa_i(K,x)\}} \prod _{i=1}^{n-1}\left( \kappa_i (K, x) -\frac{1}{R} \right)^\frac{1}{n+1} d \mu_{K\cap C}(x)\\
&&=\int_{\partial K \cap \partial C} \prod _{i=1}^{n-1} \left( \kappa_i (C, x) -\frac{1}{R} \right)^\frac{1}{n+1} d \mu_{K\cap C}(x) + \int_{\partial K \cap \partial C} \prod _{i=1}^{n-1} \left( \kappa_i (K,x) -\frac{1}{R} \right)^\frac{1}{n+1} d \mu_{K\cap C}(x).	
\end{eqnarray*}	
	%This holds since for any real numbers $a$ and $b$, we have
	%$$
	%a+b=\min\{a,b\}+\max\{a,b\}.
	%$$		
	%\vskip 2mm
%\end{proof}
\vskip 3mm
\noindent
{\bf Proof of (iii)}
We have 
\begin{eqnarray*}
as^R(K) \leq as(K) \leq n \, \vol_n(B^n_2)^\frac{2}{n+1} \vol_n(K)^\frac{n-1}{n+1}.
\end{eqnarray*}
Equality holds in the first inequality iff $R=\infty$. The second inequality is the affine isoperimetric inequality 
and equality there hold iff $K$ is an ellipsoid.
\vskip 3mm
\noindent
{\bf Proof of (iv)}
We put $f(t) = t^\frac{1}{n+1}$.  Then the proof of upper semiconituity follows immediately from \cite{Ludwig} with $f(\prod_{i=1}^{n-1} (\kappa_i(K,x) - \frac{1}{R}))$.

\end{document}